\documentclass[aps, prl, twocolumn, groupedaddress, floatfix,10pt]{revtex4-1}

\usepackage{amsmath, amssymb, amsthm}
\usepackage{graphicx}
\usepackage[varg]{txfonts} 
\usepackage{bm}
\begin{document}
\title{Functional Asynchronous Networks: Factorization of Dynamics and Function}
\author{Christian Bick${}^1$ and Michael Field${}^2$}

\affiliation{${}^1$Department of Mathematics, University of Exeter, Exeter EX4 4QF, UK\\
${}^2$Department of Mathematics, Imperial College, London SW7 2AZ, UK}
\begin{abstract}
In this note we describe the theory of
functional asynchronous networks and one of the main results, the Modularization of Dynamics Theorem,
which for a large class of functional asynchronous networks gives a factorization of dynamics in terms 
of constituent subnetworks.  For these networks we can give a complete description of the network function in terms of the function of the events comprising the network and thereby answer a question originally raised by Alon in the context of biological networks.
\end{abstract}
\maketitle

\renewcommand{\theequation}{\thesection.\arabic{equation}}
\newcommand{\NS}{N_0}
\newcommand{\conc}{\diamond}
\newcommand{\XX}{\mathbf{X}}
\newcommand{\TT}{\mathbf{T}}
\newcommand{\xx}{{\mathbf{x}}}
\newcommand{\ECS}{{\emptyset}}
\newcommand{\mc}[1]{\mathcal{#1}}
\newcommand{\is}[1]{\mathbf{#1}}
\newcommand{\iz}[1]{\is{#1}^{\bullet}}
\newcommand{\F}{\mathcal{F}}
\newcommand{\In}{\mathbb{I}}
\newcommand{\Tr}{\mathbb{F}}
\newcommand{\Net}{\mathfrak{N}}

\newcommand{\A}{\mathcal{A}}
\newcommand{\Ee}{\mathcal{E}}
\newcommand{\f}{\mathbf{f}}
\setcounter{secnumdepth}{2}
\newcommand{\dd}{\mbox{$\;|\;$}}
\newcommand{\Ref}[1]{(\ref{#1})}
\newcommand{\real}{\mbox{$\mathbb{R}$}}
\newcommand{\intg}{\mbox{$\mathbb{Z}$}}
\newcommand{\pint}{\mbox{$\mathbb{N}$}}
\newcommand{\arr}{\mbox{$\rightarrow$}}
\newcommand{\To}{\mbox{$\rightarrow$}}
\newcommand{\Mb}{\mathbf{M}}

\newcommand{\Nn}{\mathcal{N}}
\renewcommand{\aa}{\alpha}
\newcommand{\FAN}{\mathbf{N}}

\newtheorem{theorem}{Theorem}
\theoremstyle{remark}
\newtheorem{remark}{Remark}
\newtheorem{example}{Example}

\section{Introduction}
\label{intro}

Kastan \& Alon~\cite{KU} identify and describe the configurations of
relatively simple and small subnetworks that occur more frequently in biological networks than would be the case if the network
were random. 
They refer to these subnetworks as \emph{network motifs}.
Later, in his 2007 book on  systems biology~\cite{Alon}, Alon makes the following comment
\begin{quote}
\emph{``Ideally, we would like to understand the dynamics of the entire network based on the dynamics of the individual building blocks.''} (Alon~\cite[page 27]{Alon}.)
\end{quote}

The underlying premise behind this comment is that 
a modular, or engineering, approach to network dynamics is feasible. 
Identify building blocks,  connect together to form networks and then describe
dynamical properties of the resulting network in terms of the dynamics of its components.
This approach works well in linear systems theory, where a superposition principle holds, or
in, for example, the study of synchronization in weakly coupled systems of nonlinear approximately 
identical oscillators where network dynamics can be closely related
to the dynamics of individual nodes (oscillators).

However, from the perspective of a mathematician or physicist, Alon's comment seems unhelpful 
for the study of dynamics of heterogenous networks modelled by systems of nonlinear ODEs.
This is a well-known issue in complex systems~\cite{LLK}.  Yet 
engineers do couple gadgets together to make more complex systems and so it is natural
to ask if it is possible to reconcile these viewpoints.

In this note, we outline some of the basic ideas involved in an
approach to network dynamics based on what we call \emph{asynchronous networks}. We allow for features seen in networks from technology, 
engineering, and biology (for example, neuroscience or gene transcription networks). Network dynamics can involve a mix of 
distributed and decentralized control, adaptivity, event driven dynamics, switching, varying network topology and hybrid dynamics.
Typically network dynamics will be piecewise smooth, time-irreversible, nodes may stop and later restart,
and there may be no intrinsic global time.  
Significantly, many of these networks have a function: transportation networks bring people and goods from one 
point to another, neural networks may perform pattern recognition or computation, etc.
Our goal is to address Alon's comment in the context of functional asynchronous networks. Specifically, we describe 
a factorization of dynamics theorem where it is possible to describe the function of a large class of functional asynchronous networks in terms 
of the function of constituent subnetworks. As this article is only intended to be an introduction, we work always 
with the simplest examples and models and omit technical details. We refer the reader to~\cite{BF1,BF2} for greater detail, generality and proofs.

\section{Examples and properties of asynchronous networks}

We briefly describe some characteristic examples of asynchronous networks (we refer to \cite[\S 3]{BF1}
for more detail).
\subsection{Threaded \& parallel computation}
In threaded or parallel computation, computation is broken into
blocks or `threads' which are then computed \emph{independently} of each other at a speed that depends on the clock rates of the
individual processors.  As the computation proceeds, 
threads may need to exchange information with other threads. This process
involves stopping and synchronizing (updating) the thread states. Each thread may have to stop and wait for other threads to complete their
computations before it can continue with its own computation.
Each thread has its own clock (determined by its associated processor).  If threads run on separate processors, threads will be unaware of
the clock times of other threads except during the stopping and synchronization events. 

\subsection{Production and transport networks}
We give a simple detailed example of a transport network in section~\ref{sec:trans}. In production networks, parts, compounds, etc. are repeatedly built and combined as 
part of a production process leading to the desired item (for example, a car or protein). In particular, there is variation in both connection structure, typically intermittent,
and in the set of nodes (nodes may be combined or decomposed). 

\subsection{Power grid models} 
A power grid consists of a network of various types of generators and loads connected by transmission lines. 
A microgrid is a local network, capable of existing independently of the main power grid, and
typically powered by renewable energy sources (for example, solar or wind power). Critical questions involve the stability of the power
grid in case of loss of a transmission line or generator (variation in network structure), and control and stability issues related to
combining and separation (islanding) of a large set of microgrids from the main power grid. 

\subsection{Thresholds, spiking models and adaptation}  
Many mathematical models from engineering and biology incorporate thresholds. For networks,
when a node attains a threshold, there are often changes (addition, deletion, weights) in connections to another nodes.
For networks of neurons, reaching a threshold can result in a neuron firing (spiking) and short term connections to
other neurons (for transmission of the spike).  For learning mechanisms, such as Spike-Timing Dependent Plasticity (STDP)~\cite{Ge}
relative timings (the order of firing) are crucial~\cite{GK,CD,MDG} and so each neuron, or connection between a pair of neurons, comes
with a `local clock' that governs the adaptation in STDP.  In general, networks with thresholds, spiking and adaptation
provide characteristic examples of asynchronous networks where dynamics is piecewise smooth and hybrid -- a mix of
continuous and discrete dynamics.

\subsection{Properties of asynchronous networks}
We summarize some of the characteristic features of  asynchronous networks as revealed in the examples above.
\begin{enumerate}
\item Variable network structure and dependencies between nodes. Changes depend on the state of the system
or are given by a stochastic process.
\item  Synchronization events associated with stopping or waiting states of nodes.
\item  Order of events may depend on the initialization of the system or stochastic effects.
\item  Dynamics is only piecewise smooth and there may be a mix of continuous and discrete dynamics.
\item  Aspects involving function, adaptation and control.
\item  Evolution only defined for forward time -- systems are not time reversible.
\end{enumerate}

\section{Asynchronous networks: formalism}
\label{sec:async}

We use the notational conventions that $\is{k} = \{1,\ldots,k\}$ and $\iz{k} = \is{k}\cup \{0\}$, $k \in \pint$. 
Let $\real_+ = \{x \in \real \dd x > 0\}$.

Assume given a network with $k$ nodes, $N_1,\ldots,N_k$. Suppose that $N_i$ has associated phase space $M_i$, $i \in \is{k}$.
Set $\is{M} = \prod_{i \in\is{k}} M_i$ -- the network phase space.  A vector field $\f$ on  $\mathbf{M}$ is
a \emph{network vector field}.

Stopping, waiting, and synchronization are a feature of
asynchronous networks. If nodes are stopped or partially
stopped, node dynamics will be constrained to subsets of node phase
space. We codify this situation by introducing a \emph{constraining node}
$\NS$ that, when connected to $N_i$, implies that dynamics on $N_i$
is constrained. Set $\Nn = \{\NS,\ldots,N_k\}$. 

\subsection{Connection structures; admissible vector fields}

Interactions between distinct nodes in the network are given by the network graph.
Connections $N_j \arr N_i$ encode \emph{dependencies} if $i, j \in \mathbf{k}$,
and \emph{constraints} if $j = 0$, $i \in \is{k}$.

A \emph{connection structure} $\aa$ is a directed network graph on the nodes $\Nn$ such that for all $i \in \is{k}$,
$j \in \iz{k}$, $i \ne j$, there is at most one
directed connection $N_j \arr N_i$.  We do not allow self-loops or connections to the constraining node $\NS$.

An \emph{$\aa$-admissible} vector field $\f^\aa$
is a network vector field with dependencies given by $\aa$.
If $N_j \arr N_i \notin \alpha$, $i\ne j$, then
$\f^\alpha$ does not depend on $x_j\in M_j$ (and conversely, see \cite[\S 2]{BF1}). 

\begin{figure}[h]
\resizebox{0.99\columnwidth}{!}{
\includegraphics{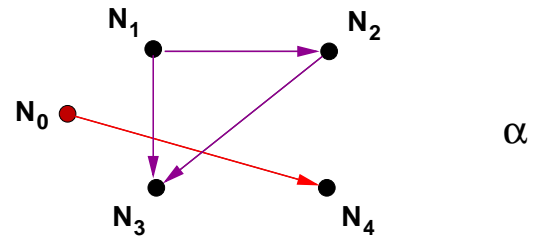} }
\caption{Construction structure on $\{\NS,N_1,N_2,N_3,N_4\}$.}
\label{fig:1}       
\end{figure}
Referring to figure~\ref{fig:1}, suppose $\f^\aa=(f_1^\aa,\ldots,f_4^\aa)$ is $\aa$-admissible. 
For  $\XX = (\xx_1,\xx_2,\xx_3,\xx_4) \in \Mb$, we have
\[
f_1^\alpha(\XX)  =  f_1(\xx_1),\; f_2^\alpha(\XX)  =  f_2(\xx_2;\xx_1)\;
\]
\[
\hspace*{-0.1in}f_3^\alpha(\XX)  =  f_3(\xx_3;\xx_1,\xx_2),\; 
f_4^\alpha(\XX)  = 0,
\]
where here we have assumed that the constraint $\NS\arr N_4$ results in $N_4$ being
stopped. 

A \emph{generalized connection structure}
$\A$ is a (nonempty) set of connection structures on $\Nn$.

An $\A$-structure $\F$ is a set
$\F = \{\f^\aa \dd \aa \in \A\}$ of network vector fields such that
each $\f^\aa\in \F$ is $\aa$-admissible.

\subsection{The event map and asynchronous networks}

Suppose that $\A$ is a generalized connection structure and $\F$ is an $\A$-structure.

Interactions between nodes in asynchronous networks
may be state dependent or change over time (stochastically).
Here we only consider state dependence.

We handle interactions and constraints using an \emph{event map} $\mathcal{E}: \Mb \arr \A$.

The quadruple $\Net = (\Nn,\A,\F,\Ee)$
defines an \emph{asynchronous network}.
Dynamics on $\Net$ is
given by the state dependent network vector field $\F$ defined by
\begin{equation}
\label{eq:1}
\is{F}(\XX) = \f^{\Ee(\XX)}(\XX), \; \XX \in \Mb.
\end{equation}
\begin{remark}
We refer to~\cite[\S 4.7]{BF1} for the definition of an integral curve for \Ref{eq:1}
and note that the definition we use is different from that used in Filippov systems~\cite{Fil,BBCK}.
Under simple conditions on the event map~\cite[\S\S 4.7, 4.8]{BF1}, it can be shown that if $\Mb$ is compact then for each $\XX \in \Mb$
there is a unique piecewise smooth integral curve $\Phi_\XX:[0,\infty)\arr \Mb$ with initial condition $\XX$ and corresponding
semiflow $\Phi:\Mb\times\real_+ \arr \Mb$. In general, $\Phi$ will not be continuous as a function of $\XX\in\Mb$. 
\end{remark}

\section{A simple transport example}
\label{sec:trans}
\begin{figure}[h]
\resizebox{0.99\columnwidth}{!}{
\includegraphics{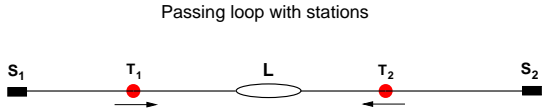} }
\caption{A single track railway line with a passing loop.}
\label{fig:2}
\end{figure}
We consider a single track railway line joining stations marked $S_1, S_2$ in figure~\ref{fig:2}.
Suppose there is a passing loop at $L$. Trains, marked $T_1$ and $T_2$ in figure~\ref{fig:2}, start
from stations $S_1$ and $S_2$ respectively and proceed towards the opposite station. There is no
central control or communication between the train drivers except when both trains are in the passing loop.
We further assume that both drivers are running a nonlinear oscillator.  
When a train enters the passing loop it stops. When both trains are in the passing loop, the 
drivers cross couple their nonlinear oscillators.
In order for a train to leave the passing loop, two conditions must be satisfied.
\begin{enumerate}
\item Both trains must be in the passing loop.
\item The two nonlinear oscillators must be phase synchronized to within $\varepsilon$ where
$0 < \varepsilon \ll \pi$.
\end{enumerate}
We model this setup using an asynchronous network with two nodes -- corresponding to the two trains.
The phase space for train $T_i$ will be $M_i = [-a,b]\times \mathbb{T}$, $i\in\is{2}$ ($\mathbb{T}=\real/2\pi \intg$),
where the interval $[-a,b]$ models the line joining $S_1$ to $S_2$, $S_1$ has coordinate $-a < 0$, $S_2$ has coordinate $b > 0$,
The passing loop $L$ will be at $x = 0$, and $\mathbb{T}$ will be the phase space for the nonlinear oscillator.

Assume train motion given by $V_1, V_2$ where $V_1(x) > 0 > V_2(x)$ $x \in [-a,b]$.   Note that neither $V_1$ or $V_2$ can be zero anywhere on $[-a,b]$ otherwise the trains
will never both reach their destination stations in finite time. 

We define four connection structures.

\begin{eqnarray*}
\alpha_i&=&N_0\arr N_i,\; i\in\is{2} \;(\mbox{$T_i$ stopped})\\
\beta &= & N_0\arr N_1\leftrightarrow N_2\leftarrow N_0(\mbox{stopped \& cross coupled})\\ 
\ECS & = &  \mbox{Empty connection structure} 
\end{eqnarray*}
Let $\A=\{\alpha_1,\alpha_2,\beta,\ECS\}$ be the associated generalized connection structure.

Model the uncoupled oscillator dynamics for train $T_i$ by $\dot{\theta}_i = \omega$, where $\omega > 0$, and the coupled dynamics
by the Kuramoto phase oscillator system
\begin{eqnarray*}
\dot{\theta}_1 & = & \omega + \sin(\theta_2-\theta_1)\\
\dot{\theta}_2 & = & \omega + \sin(\theta_1-\theta_2)
\end{eqnarray*}

It remains to define the admissible vector fields and event map that give the required dynamics for this example.
As admissible vector fields (on $([-a,b] \times \mathbb{T})^2$) we take
\begin{eqnarray*}
\is{f}^\ECS\,  &=& ((V_1,g),(V_2,g)),\\
\is{f}^{\alpha_1} & = & ((0,g),(V_2,g)),\;(T_1\;\mbox{stopped}, \; T_2\;\mbox{running})\\
\is{f}^{\alpha_2} & = & ((V_1,g),(0,g)),\;(T_2\;\mbox{stopped}, \; T_1\;\mbox{running})\\
\is{f}^\beta\, & =& ((0,G_1), (0,G_2)), \; (T_1,T_2\;\mbox{stopped \& cross coupled})
\end{eqnarray*}
We define the event map by
\begin{eqnarray*}
\Ee(\XX,\bm{\theta}) &=& \aa_1,\; x_1 = 0, x_2 > 0\\
& = & \aa_2, \; x_1<0, x_2=0\\
& = & \beta,\; x_1 = x_2=0,\; |\theta_1-\theta_2| > \varepsilon \\
& = & \ECS, \;\mbox{otherwise}
\end{eqnarray*}
Finally, dynamics are given by the network vector field

\begin{equation}
\label{eq:2}
F(\XX,\bm{\theta}) = \is{f}^{\Ee(\XX,\bm{\theta})}(\XX,\bm{\theta}).
\end{equation}
We leave it as an easy exercise for the reader to check that 
if $T_i$ is initialized at $(x_i(0),\theta_i(0))\in [-a,b]\times \mathbb{T}$, $i \in \is{2}$,
where $x_1(0) \le 0 \le x_2(0)$, then \Ref{eq:2}
has a well defined integral curve $\varphi: \real_+ \arr ([-a,b]\times \mathbb{T})^2$,  with specified initial condition,
that gives the correct dynamics for the passing loop problem.

\begin{remark}
The passing loop gives an example of a simple \emph{functional} asynchronous network. The function is 
for the trains to go from one station to the opposite station in finite time. Observe that for this example there is the possibility of a
\emph{dynamical deadlock}: if the trains start at the same time and if $\theta_1(0) = \theta_2(0) + \pi$, then the 
coupled phase oscillators will never phase synchronize -- $\theta_1(t) = \theta_2(t) + \pi$ for all $t \in \real_+$ -- and so the trains
will never exit the passing loop. We refer to~\cite[\S\S 2,3]{BF2} for more details on deadlocks in functional asynchronous networks.
\end{remark}
\section{Functional asynchronous networks}
We follow the notational conventions of section~\ref{sec:async} and let 
$\Net = (\mathcal{N},\mathcal{A},\mathcal{F},\mathcal{E})$ denote an asynchronous network. We assume
that $\Net$ has associated semiflow 
\[
\Phi = (\Phi_1,\ldots,\Phi_k): \Mb\times \real_+\arr \Mb.
\]
Suppose that we are given 
\emph{initialization} and \emph{termination} sets $\In,\Tr \subset \Mb$ where
\[
\In  =  \prod_{i \in \is{k}}\In_i,\quad \Tr  =  \prod_{i \in \is{k}}\Tr_i,
\]
Typically, $\In_i,\Tr_i\subset M_i$ will be closed disjoint hypersurfaces that separate $M_i$ into three connected components, $i \in \is{k}$. That is,
$M_i = M_i^- \cup M_i^0 \cup M_i^+$ where $M_i^- \cap M_i^+ = \emptyset$ and
\[
M_i^- \cap M_i^0 = \partial M_i^- = \In_i,\quad M_i^0 \cap M_i^+ = \partial M_i^+ = \Tr_i.
\]
We call $\FAN = (\Net,\In,\Tr)$ a
\emph{functional asynchronous network}. The network function is getting from $\In$ to $\Tr$ and is expressed by the transition
and timing functions
\[
G_0: \mc{D}\subset \In\arr \Tr,\quad \is{S}: \mc{D}\subset \In \arr \real_+^k.
\]
That is, if $\XX \in \mc{D}$, then for all $i \in \is{k}$ there exists $S_i \in \real_+^k$ such that
\begin{eqnarray*}
\Phi_i(\XX,S_i)&\in&\Tr_i, \quad \Phi_i(\XX,t) \notin \Tr_i, t < S_i,\\
\is{S}(\XX) & = & (S_1,\ldots,S_k),
\end{eqnarray*}
\begin{example}
For the passing loop example discussed in the previous section, we take
$\In_1 = M_1^- = M_2^+ = \Tr_2 = \{-a\}\times \mathbb{T}$, $\Tr_1= M_1^+ = M_2^+ = \In_2 = \{b\}\times\mathbb{T}$.
In this case, there is the implicit assumption that trains stop when they reach their termination set.
Equally well, we could take $M_i = \real\times\mathbb{T}$ so that $M_1^- = (-\infty,a]\times\mathbb{T}$ etc. (see also~\cite[\S3]{BF2}).
Finally, observe that $\mc{D} = \{((-a,\theta_1),(b,\theta_2))\dd |\theta_1-\theta_2| \ne \pi\}$.
\end{example}

More generally, we allow for general initialization times and define generalized transition and timing functions
\[
G: \widehat{\mc{D}}\subset \In\times \real_+^k\arr \Tr,\quad \widehat{\is{S}}: \widehat{\mc{D}}\subset \In\times \real_+^k \arr \real_+^k.
\]
We refer to \cite[\S 3.4]{BF2} for details. For our main result, it is required that the network has a
generalized transition and timing functions with $\widehat{\mc{D}} = \In\times \real_+^k$.

\subsection{Functional networks built from events}
\begin{figure}[h]
\resizebox{0.99\columnwidth}{!}{
\includegraphics{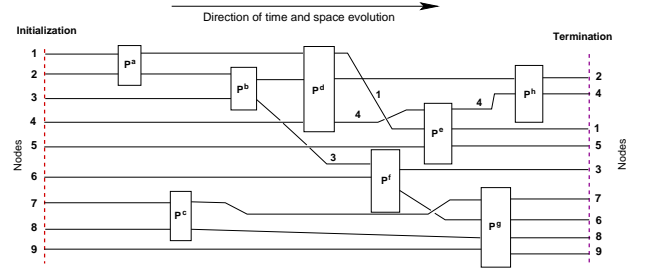} }
\caption{A spatiotemporal decomposition of a functional asynchronous network}
\label{fig:3}
\end{figure}

In figure~\ref{fig:3} we show a nine node functional asynchronous network that is built from the eight ``events'' $\is{P}^a,\ldots,\is{P}^h$.

The initialization and termination sets are indicated on the left and right sides of the figure respectively.
The events signify regions of phase space where there can be (state dependent) interaction between nodes. For example, the event labelled
$\is{P}^g$ involves interaction between nodes $N_6$, $N_7$, $N_8$, and $N_9$. Observe that there is only a partially ordered temporal structure 
on the events. Thus, the event $\is{P}^g$ must occur after $\is{P}^f$ but can occur before or after event $\is{P}^h$.

\subsection{Building blocks}
\begin{figure}[h]
\resizebox{0.99\columnwidth}{!}{
\includegraphics{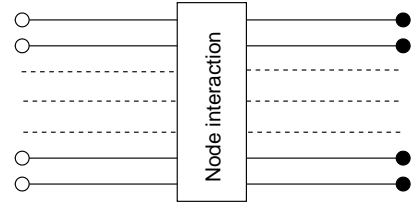} }
\caption{Dynamical/functional module}
\label{fig:4}
\end{figure}

In figure~\ref{fig:4} we represent a basic building block with the same number of inputs and outputs.

The initialization sets are represented by the symbols $\circ$, termination sets by $\bullet$. Interaction between nodes occurs only in the event
region denoted by the rectangle. Outside of the event region, nodes evolve independently. 
More generally, we can allow for different number of inputs and outputs: nodes may merge or split.

Our immediate aim is describe some basic operations that we can define on functional asynchronous networks that allow us enable us to find
a (maximal) decomposition of a 
functional asynchronous network into the form shown in figure~\ref{fig:3}.
\subsection{Operations on functional asynchronous networks}

If $\FAN^a = (\Net^a,\In^a,\Tr^a)$, $a\in \is{q}$, are functional asynchronous networks with distinct node sets ($\Nn^a \cap\Nn^b \subset \{\NS\}$, $a\ne b \in \is{q}$),
define the \emph{product} $\prod_{a\in\is{q}}\FAN^a$ to be the functional asynchronous network $\FAN = (\Net,\In,\Tr)$, where
\[
\In = \prod_{a\in\is{q}} \In^a,\quad \Tr = \prod_{a\in\is{q}} \Tr^a
\]
and $\Net = \prod_{a\in \is{q}} \Net^a$ is defined in the obvious way to be the asynchronous network with node set $\Nn = \cup_{a\in\is{q}}\Nn^a$ 
(we refer to \cite[\S 6]{BF1} for details).

We say the $k$-node functional asynchronous network $\FAN = (\Net,\In,\Tr)$ is \emph{trivial} if $\FAN = \prod_{a \in \is{k}} \FAN^a$ where
each $\FAN^a$ has exactly one node $N_a$. In particular, if $\FAN$ is trivial there are no interactions between nodes and no constraints.

Next let $\Net^a = (\Net^a,\In,\Tr)$, $a\in \is{q}$, be a family of functional asynchronous networks with common initialization set,  termination set and node set
$\Nn = \{\NS,N_1,\ldots,N_k\}$. Suppose that for each
$a \in \is{q}$, there exists $\Sigma(a) \subset \is{k}$ such that 
\begin{enumerate}
\item $\Net^a = \Net^a_1 \times \Net^a_2$ where $\Net^a_1$ has node set $\Sigma(a)$ and $\Net^a_2$ is trivial.
\item If $a\ne b$, $\Sigma(a) \cap \Sigma(b) = \emptyset$.
\end{enumerate}
We define the \emph{amalgamation} $\Net = \sqcup_{a\in\is{q}} \Net^a$ to be the functional asynchronous network $(\prod_{a\in \is{q}}  \Net^a_1)  \times \Net_2$,
where $\Net_2$ is the trivial network defined as the product of the common trivial factors in $\Net^a_2$, $a \in \is{q}$. Thus the node set
of $\Net_2$ will be $\is{k} \smallsetminus \cup_{a\in\is{q}} \Sigma(a)$.
\begin{figure}[h]
\resizebox{0.99\columnwidth}{!}{
\includegraphics{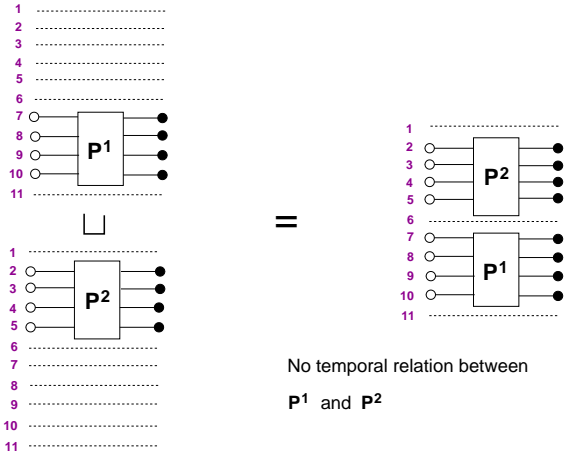} }
\caption{Amalgamating two functional asynchronous networks.}
\label{fig:5}
\end{figure}

Referring to figure~\ref{fig:5}, we have $\Sigma(1) = \{7,8,9,10\}$ and $\Sigma(2) = \{2,3,4,5\}$.
The amalgamation $\is{P} = \is{P}^1 \sqcup \is{P}^2$ is trivial when restricted to nodes $\{N_1,N_6,N_{11}\}$.

\begin{figure}[h]
\resizebox{0.99\columnwidth}{!}{
\includegraphics{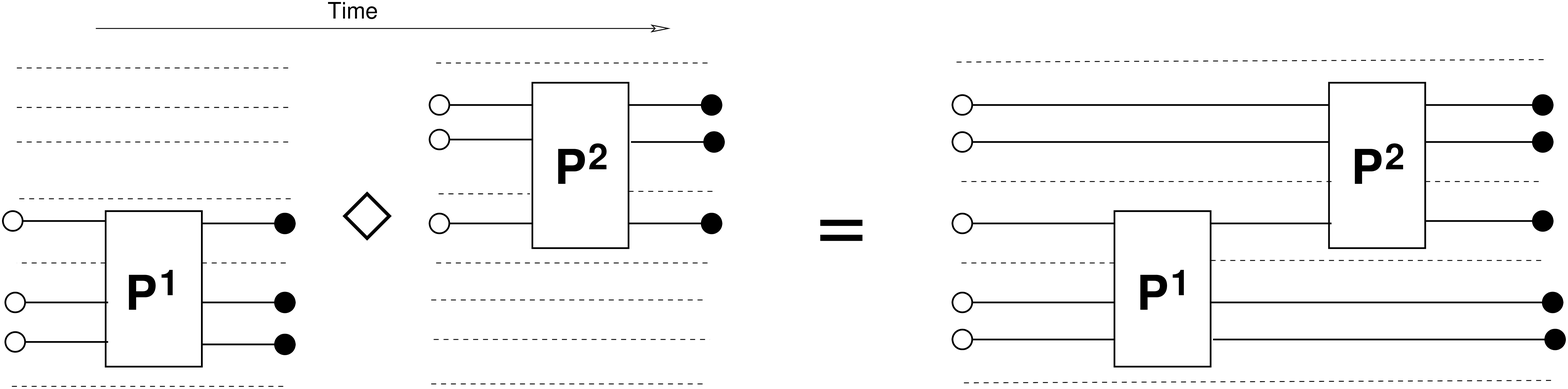} }
\caption{Concatenating two functional asynchronous networks.}
\label{fig:6}
\end{figure}

Finally, we outline the operation of concatenation, referring the reader to~\cite[\S4]{BF2} for the details (most) we omit.
Suppose that $\Net^a = (\Net^a,\In^a,\Tr^a)$, $a\in \is{2}$, are functional asynchronous networks with common node set. 
Assume that $ \Tr^1 = \In^2$. The \emph{concatenation} $\FAN = (\Net,\In,\Tr) = \FAN^2 \conc \FAN^1$ will be a temporal merging $\FAN^1, \FAN^2$.
We define 
\begin{enumerate}
\item $\In = \In^1$, $\Tr = \Tr^2$.
\item $\A = \{\aa_1 \vee \aa_2 \dd \exists \XX\in\Mb, \alpha_1 = \Ee^1(\XX), \aa_2 = \Ee^2(\XX)\}$,
\end{enumerate}
where $\vee$ denote the join of the graphs. The definition of the set of admissible vector fields $\F$
for $\Net$ is trickier and requires additional conditions on $\FAN^1, \FAN^2$ -- we refer to ~\cite{BF2} for details.  
We define the event map by $ \Ee(\XX) = \Ee^1(\XX) \vee \Ee^2(\XX), \XX \in \Mb$. We refer to figure~\ref{fig:6}
for the operation of concatenation.

The concatenation $\FAN^2 \conc \FAN^1$ has the important property that if  $\FAN^a$ has generalized transition and timing functions 
$G^a: \In^a \times \real_+^k \arr \Tr^a$, $\widehat{\is{S}}^a: \In^a \times \real_+^k\arr \real_+^k$, $a \in \is{2}$, then 
$\FAN^2 \conc \FAN^1$ has generalized transition function $G$ given by
$G(\XX,\TT) = G^2(G^1(\XX,\TT), \widehat{\is{S}}^1(\XX,\TT))$~\cite[Corollary 4.15]{BF2}.

\begin{remark}
We have deliberately avoided listing the detailed properties required of functional asynchronous 
networks in order to define amalgamations and concatenations. Briefly, apart from requiring
the existence of generalized transition and timing functions, we require (1) the uncoupled vector
vectors defining intrinsic dynamics of a node $N_i$ to be transverse to $\In_i$, $\Tr_i$ and (2)
a local product structure on the network.  We refer to~\cite[\S 3]{BF2} for the details.
\end{remark}
\section{Modularization of dynamics and function}
A functional asynchronous network is \emph{primitive} if it cannot be written as a nontrivial amalgamation or concatenation.
\begin{theorem}
\label{fan}
Under general conditions, a functional asynchronous network $\FAN$ has a
unique (up to rearrangements) decomposition
\[
\FAN = \FAN^q \conc \ldots \conc \FAN^1,
\]
where
$ \FAN^j = \FAN^{j,1} \sqcup \ldots \sqcup \FAN^{j,q(j)}$, $j \in \is{q}$, and 
and each $\FAN^{j,\ell}$ is primitive.

The generalized transition function $G$ for $\FAN$ can be expressed in terms of the
generalized transition and timing functions $G^j, \widehat{\is{S}}_j$ of $\FAN^j$ (or $G^{j,\ell}$ for $\FAN^{j,\ell}$) by:
\begin{eqnarray*}
G(\XX,\TT)&=&G^q(\ldots G^2(G^1(\XX,\TT),\widehat{\is{S}}^1(\XX,\TT))\ldots),\\ 
\widehat{\is{S}}(\XX,\TT) & = & \widehat{\is{S}}^q(\ldots \widehat{\is{S}}^2(G^1(\XX,\TT), \widehat{\is{S}}^1(\XX,\TT))\ldots),\\
G^j&=&G^{j, 1} \times \ldots \times G^{j,q(j)}, \; j\in \is{q}.
\end{eqnarray*}
\end{theorem}
\begin{example}                                 
Consider the network shown in figure~\ref{fig:3} and assume that each event $\is{P}^j$. $j \in \{a,\ldots,h\}$ is primitive.
A decomposition satisfying the requirements of theorem~\ref{fan} is indicated in figure~\ref{fig:7} -- the dashed lines 
indicate the initialization and termination sets for the subnetworks.
\begin{figure}[h]
\resizebox{0.99\columnwidth}{!}{
\includegraphics{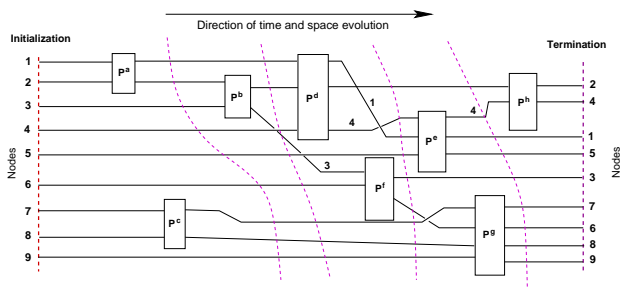} }
\caption{Factorization of network of figure~\ref{fig:3}.}
\label{fig:7}
\end{figure}
The factorization for the network is
\[
\FAN = \is{P}^h\conc (\is{P}^e\sqcup\is{P}^g)\conc(\is{P}^d \sqcup \is{P}^f)\conc\is{P}^b\conc(\is{P}^a \sqcup \is{P}^c).
\]
This factorization corresponds to maximizing from the left hand side.
However, if we maximize from the right we obtain the factorization
\[
\FAN = (\is{P}^h\sqcup\is{P}^g) \conc (\is{P}^c \sqcup\is{P}^e\sqcup\is{P}^f)\conc\is{P}^d \conc\is{P}^b\conc\is{P}^a.
\]
In either case there is a concatenation of five networks -- that is the minimum number possible.
\end{example}

Theorem~\ref{fan} allows us to write the function of a network explicitly in terms of the 
transition functions of the constituent subnetworks.  

Results of this type depend crucially on intermittent connection structure and nonsmooth dynamics. For example, no such
result is possible for a classical coupled network of phase oscillators.

The approach works because we have adopted an engineer's viewpoint: we emphasise function rather than dynamics. Indeed,
we are indifferent to the specific dynamics occurring \emph{between} the initialization and termination sets. Of course, both the timing and
transition functions provide the key information about network function. 

\section{Concluding comments}
\begin{enumerate}

\item Theorem~\ref{fan} is a prototypical theorem providing proof of concept.
The conditions for the theorem can be significantly weakened from those required in~\cite{BF2}.

\item The theorem yields maximal feedforward structure on a functional asynchronous network
(note that individual events may have feedback loops).

\item The result suggests the utility of starting with a small functional asynchronous network; understanding the structure in depth and then 
then evolving to optimize function (for example by adding feedback).

\item There are many as yet unexplored issues such as bifurcation, hidden deadlocks, and the effects of noise.

\item There is the problem of how far one can determine internal structure on the basis of input/output time series data.

\end{enumerate}


\begin{thebibliography}{99}
\bibitem{Alon} U Alon. \emph{An Introduction to Systems Biology. Design Principles of Biological Circuits} (Chapman \& Hall/CRC, Boca Raton, 2007).
\bibitem{BBCK} M di Bernardo, C J Budd, A R Champneys, \& P Kowalczyk. \emph{Piecewise-smooth Dynamical Systems} (Springer, Applied Mathematical Science 163, London, 2008).
\bibitem{BF1} C Bick \& M J Field. `Asynchronous networks and event driven dynamics', (2016), preprint.
\bibitem{BF2} C Bick \& M J Field. `Asynchronous networks: Modularization of Dynamics Theorem', (2016), preprint.
\bibitem{CD} N Caporale \& Y Dan. `Spike timing dependent plasticity: a Hebbian learning rule' \emph{Annu. Rev. Neurosci.} {\bf 31} (2008), 25--36.
\bibitem{Fil} A F Filippov. \emph{Differential Equations with Discontinuous Righthand Sides} (Kluwer Academic Publishers, 1988).
\bibitem{Ge} W Gerstner, R Kempter, L J van Hemmen \& H Wagner. `A neuronal learning rule for sub-millisecond temporal coding', \emph{Nature} {\bf 383}
(1996), 76--78.
\bibitem{GK} W Gerstner \& W Kistler. \emph{Spiking Neuron Models} (Cambridge University Press, 2002).
\bibitem{KU} N Kashtan \& U Alon. `Spontaneous evolution of modularity and network motifs',
\emph{PNAS} {\bf 102} (39) (2005), 13773--13778.
\bibitem{LLK} J Ladyman, J Lambert, and K Wiesner. `What is a complex system?', \emph{Eur. J. for Phil. of Sc.} {\bf 3}(1) (2013), 33--67.
\bibitem{MDG} A Morrison, M Diesmann \& W Gerstner. `Phenomenological models of synaptic plasticity based on spike timing',
\emph{Biol. Cyber.} {\bf 98} (2008), 459--478.
\end{thebibliography}
\end{document}